\documentclass[reqno]{amsart}

\usepackage{latexsym,amssymb,amsthm,amsmath}
\usepackage{graphicx,url,xcolor}

\usepackage{xypic}
\usepackage[arc,all]{xy}

\theoremstyle{plain}

\newtheorem{theorem}{Theorem}

\newtheorem{proposition}{Proposition}
\newtheorem{corollary}{Corollary}

\newtheorem{remark}{Remark}

\begin{document}

\title[Semi-invariant $\xi ^{\bot}$-submanifolds of G. Q. S. manifolds]
{Semi-invariant $\xi ^{\bot}$-submanifolds of generalized quasi-Sasakian manifolds}
\author[C.~C\u alin]{Constantin C\u alin}
\address[C.~C\u alin]{Technical University Gh. Asachi,
Department of Mathematics, Iasi, 700049
Romania}
\email{c0nstc (at) yahoo.com}

\author[M.~Cr\^a\c smareanu]{Mircea Cr\^a\c smareanu}
\author[M.~I.~Munteanu]{Marian Ioan Munteanu}
\thanks{The third author was supported by Grant PN-II ID 398/2007-2010 (Romania)}
\address[M.~Cr\^a\c smareanu, M.~I.~Munteanu]{'Al.I.Cuza' University of Iasi\\ Bd. Carol I, nr. 11\\
Iasi, 700506 \\ Romania,\quad
\url{http://www.math.uaic.ro/~mcrasm}\ ,
\
\url{http://www.math.uaic.ro/~munteanu}
}
\email{mcrasm (at) uaic.ro, marian.ioan.munteanu (at) gmail.com}

\author[V.~Saltarelli]{Vincenzo Saltarelli}

\address[V.~Saltarelli]{Department of Mathematics,
University of Study of Bari, Via E. Orabona 4, 70125 - Bari, Italy  }
\email{saltarelli (at) dm.uniba.it}

\date{\today}

\dedicatory{
Dedicated to the memory of Prof. Stere Ianu\c s $(1939-2010)$
}
\begin{abstract}
A structure on an almost contact metric manifold is defined as a
generalization of well-known cases:
Sasakian, quasi-Sasakian, Kenmotsu and cosymplectic.
Then we consider a semi-invariant $\xi^{\bot}$-submanifold
of a manifold endowed with such a structure and two topics are studied:
the integrability of distributions defined by this submanifold and
characterizations for the totally umbilical case. In particular
we recover results of Kenmotsu \cite{5}, Eum \cite{7} and Papaghiuc \cite{9}.
\end{abstract}



\maketitle

\section{Preliminaries and basic formulae}

An interesting topic in the differential geometry is the theory of submanifolds in
spaces endowed with additional structures. In 1978, A.~Bejancu (in \cite{2}) studied CR-submanifolds
in K\"ahler manifolds. Starting from it, several papers have been appeared in this field.
Let us mention only few of them: a series of papers of B.Y.~Chen (e.g. \cite{6}), of A.~Bejancu and N.~Papaghiuc (e.g. \cite{3}
in which the authors studied semi-invariant submanifolds in Sasakian manifolds).
See also \cite{Mun05}. The study was extended also to other ambient spaces, for example
A.~Bejancu in \cite{4} also studied QR-submanifolds in quaternionic manifolds and
M.~Barros in \cite{1} investigated CR-submanifolds in quaternionic manifolds. Several important
results above CR-submanifolds are being brought together in \cite{4}, \cite{6}, \cite{Man04},
\cite{Mun05}, \cite{Orn06} and the corresponding references.
The purpose of the present paper is to investigate the
semi-invariant $\xi^{\bot}$-submani\-folds in a generalized Quasi-Sasakian manifold.

Let $\widetilde{M}$ be a real $(2n+1)$-dimensional smooth manifold
endowed with an almost contact metric structure $(\phi,\xi ,\eta ,\tilde g)$:
$$
\left\{
\begin{array}{l}
\phi^{2}=-I+\eta \otimes \xi ,\ \eta (\xi )=1,\  \eta\circ \phi=0, \ \phi\xi =0\\
\eta (X)=\tilde g(X,\xi ),\tilde g(\phi X,Y)+\tilde g(X,\phi Y)=0
\end{array}
\right.
$$
for any vector fields $X,Y$ tangent to $\widetilde{M}$ where $I$ is the identity
on sections of the tangent bundle $T\widetilde{M}$, $\phi$ is a tensor field of
type $(1,1)$, $\eta $ is a 1-form, $\xi $ is a vector field and $\tilde g$ is a Riemannian metric
on $\widetilde{M}$.
Throughout the paper all manifolds and maps are smooth. We denote by ${\mathcal F}(\widetilde{M})$ the algebra
of the smooth functions on $\widetilde{M}$ and by $\Gamma (E)$ the
${\mathcal F}(\widetilde{M})$-module of the sections of a vector bundle $E$
over $\widetilde{M}$.

The almost contact manifold $\widetilde{M}(\phi,\xi ,\eta)$ is said to be
{\it normal} if
\[
N_{\phi}(X,Y)+2d\eta (X,Y)\xi =0
\]
where
\[
N_{\phi}(X,Y)=[\phi X,\phi Y]+\phi^{2}[X,Y]-\phi [\phi X,Y]-\phi [X,\phi Y],\quad X,Y\in \Gamma (T\widetilde{M})
\]
 is the Nijenhuis tensor field corresponding of the tensor field $\phi$.

The {\it fundamental $2$-form} $\Phi$ on $\widetilde{M}$ is defined by
$\Phi (X,Y)=\tilde g(X,\phi Y)$.

In \cite{5}, the author studied hypersurfaces of an almost
contact metric manifold $\widetilde{M}$ whose structure tensor fields
satisfy the following relation
\begin{equation}
\label{(1.1)}
(\widetilde{\nabla}_{X}\phi)Y=\tilde g(\widetilde{\nabla}_{\phi X}\xi ,Y)\xi -\eta (Y)\widetilde{\nabla }_{\phi X}\xi
\end{equation}
where $\widetilde{\nabla}$ is the Levi-Civita connection of the metric tensor $\tilde g$.
See also \cite{7,8}.
For the sake of simplicity we say that a manifold $\widetilde{M}$ endowed with an
almost contact metric structure satisfying \eqref{(1.1)} is a {\it generalized Quasi-Sasakian
manifold}, in short G.Q.S. Define a $(1,1)$ type tensor field $F$ by
\begin{equation}
\label{(1.2)}
FX=-\widetilde{\nabla}_{X}\xi .
\end{equation}

\begin{proposition}
\label{p:1.1}  If $\widetilde{M}$ is a G.Q.S manifold then any
integral curve of the structure vector field $\xi $ is a geodesic i.e.
$\widetilde{\nabla}_{\xi}\xi =0$. Moreover $d\Phi =0$ if and only if $\xi $ is a
Killing vector field.
\end{proposition}

\proof The first assertion follows immediately from \eqref{(1.1)} with $X=Y=\xi$, and taking into account
that $\eta(\widetilde\nabla_\xi\xi)=0$. Next, we deduce
$$
3d\Phi(X,Y,Z)=\tilde g\big((\widetilde{\nabla}_X\phi)Z,Y\big)+\tilde g\big((\widetilde{\nabla}_Z\phi)Y,X\big)+
\tilde g\big((\widetilde{\nabla}_Y\phi)X,Z\big)+
$$
$$
+\eta(X)\Big(\tilde g(Y,\widetilde{\nabla}_{\phi Z}\xi)+\tilde g(\phi Z,\widetilde{\nabla}_Y\xi)\Big)+
\eta(Y)\Big(\tilde g(Z,\widetilde{\nabla}_{\phi X}\xi)+\tilde g(\phi X,\widetilde{\nabla}_Z\xi)\Big)+
$$
$$
+\eta(Z)\Big(\tilde g(X,\widetilde{\nabla}_{\phi Y}\xi)+\tilde g(\phi Y,\widetilde{\nabla}_X\xi)\Big).
$$
If we suppose that $\xi $ is Killing then, from the last equation, we obtain $d\Phi =0$.

Conversely, suppose that $d\Phi =0$. Taking into account the first part of the statement,
for $X=\xi$, $\eta(Y) = \eta(Z) =0$, the last relation implies
$$
\tilde g\big(Y, \widetilde{\nabla}_{\phi Z}\xi\big) +\tilde g\big (\phi Z, \widetilde{\nabla}_Y\xi\big) =0.
$$
Finally, by replacing $Z$ with $\phi Z$ and $Y$ by $Y-\eta(Y)\xi $ we deduce that
$\xi $ is a Killing vector field.
\endproof

The next result can be obtained by direct calculation:

\begin{proposition}
\label{p:1.2}
A G.Q.S manifold $\widetilde{M}$ is normal and
\begin{equation}
\label{(1.3)}
\phi\circ F=F\circ \phi, \ F\xi =0,\
 \eta \circ F=0,\ \widetilde{\nabla}_{\xi }\phi=0.
\end{equation}
\end{proposition}

\begin{remark}
\label{r:1.3}
\rm
{\bf a)} It is easy to see that on such manifold $\widetilde{M}$ the
structure vector field $\xi $ is not necessarily a Killing vector field i.e. $\widetilde{M}$
is not necessarily a K-contact manifold.
{\bf b)} It is also interesting to pointed out that the following particular situations hold
\begin{itemize}
\item [1)] $FX= -\phi X$ then $\widetilde{M}$ is Sasakian
\item [2)] $FX =-X + \eta(X)\xi $ then $\widetilde{M}$ is Kenmotsu
\item [3)] $FX = 0$ then $\widetilde{M}$ is cosymplectic
\item [4)] if $\xi $ is a Killing vector field then $\widetilde{M}$ is a quasi-Sasakian manifold.
\end{itemize}
\end{remark}

Now, let $\widetilde{M}$ be a G.Q.S manifold and consider an $m$-dimensional submanifold $M$, isometrically immersed in $\widetilde{M}$.
Denote by $g$ the induced metric on $M$ and by $\nabla $ its Levi-Civita connection.
Let $\nabla^{\bot}$ and $h$ be the normal
connection induced by $\widetilde{\nabla}$ on the normal bundle $TM^{\bot}$ and
the second fundamental form of $M$, respectively. Then one has the direct sum decomposition
$T\widetilde{M} = TM \oplus TM^{\bot} $.
Recall the Gauss and Weingarten formulae
$$
\begin{array}{l}
{\rm (G)} \ \ \widetilde{\nabla}_{X}Y = \nabla_{X}Y + h(X,Y)  \\[2mm]
{\rm (W)} \ \widetilde{\nabla}_{X}N = - A_NX + \nabla_X^{\bot}N, \quad X, Y \in \Gamma(TM)
\end{array}
$$
where $A_N$ is the shape operator with respect to the normal section $N$
and satisfies
\[
\tilde g(h(X, Y), N) = g(A_NX, Y) \quad X, Y\in\Gamma(TM), \ \
N\in\Gamma(TM^{\bot}).
\]

The purpose of the present paper is to investigate the
semi-invariant $\xi^{\bot}$-submani\-folds in a G.Q.S manifold. More
precisely, we suppose that the structure vector field $\xi$ is orthogonal to
the submanifold $M$. According to Bejancu \cite{4} we say that $M$ is a {\it semi-invariant
$\xi^{\bot}$-submanifold} if there exist two orthogonal distributions, ${\mathcal{D}}$
and ${\mathcal{D}}^{\bot}$, in $TM$ such that:
\begin{equation}
\label{(1.4)}
 TM = {\mathcal{D}}\oplus {\mathcal{D}}^{\bot},  \ \phi {\mathcal{D}} = {\mathcal{D}}, \
\phi {\mathcal{D}}^{\bot}\subseteq TM^{\bot}
\end{equation}
where $\oplus$ denotes the orthogonal sum. If ${\mathcal{D}}^{\bot} = \{0\}$ then $M$ is
an {\it invariant $\xi^{\bot}$-submanifold}.
The normal bundle can also be decomposed as $TM^\bot=\phi{\mathcal{D}}^\bot\oplus\mu$,
where $\phi\mu\subseteq\mu$. Hence $\mu$ contains $\xi$.

\section{Integrability of distributions on a se\-mi-in\-variant $\xi^{\bot}$-submanifold}

Let $M$ be a semi-invariant $\xi^{\bot}$-submanifold of a G.Q.S
manifold $\widetilde{M}$. Denote by $P$ and $Q$ the projections of $TM$
on ${\mathcal{D}}$ and ${\mathcal{D}}^{\bot }$ respectively, namely for any $X\in \Gamma (TM)$
\begin{equation}
\label{(2.1)}
X=PX+QX.
\end{equation}
Moreover, for any $X\in \Gamma (TM)$ and $N\in \Gamma (TM^{\bot })$ we put
\begin{equation}
\label{(2.2)}
\phi X=tX+\omega X
\end{equation}
\begin{equation}
\label{(2.3)}
\phi N=BN+CN
\end{equation}
with $tX\in \Gamma ({\mathcal{D}})$, $BN\in \Gamma (TM)$ and $\omega X, CN\in \Gamma (TM^{\bot })$.
We also consider, for $X\in \Gamma (TM)$, the decomposition
\begin{equation}
\label{(2.4)}
FX=\alpha X+\beta X,\ \ \ \alpha X\in \Gamma ({\mathcal{D}}),\ \beta X\in \Gamma
(TM^{\bot }).
\end{equation}

The purpose of this section is to study the integrability of both
distributions ${\mathcal{D}}$ and ${\mathcal{D}}^{\bot }$. With this scope in mind, we state
first the following result.

\begin{proposition}
\label{p:2.1}
Let $M$ be a semi-invariant $\xi^{\bot}$-submanifold of a G.Q.S manifold $\widetilde{M}$. Then we have
\begin{equation}
\label{(2.5)}
\begin{array}{l}
a)\ (\nabla _{X}t)Y=A_{\omega Y}X+Bh(X,Y),\vspace{2mm} \\
b)\ (\nabla _{X}\omega )Y=Ch(X,Y)-h(X,tY)+g(FX,\phi Y)\xi ,\quad X,Y\in \Gamma (TM).
\end{array}
\end{equation}
\end{proposition}
\proof
The statement follows immediately from \eqref{(2.2)}--\eqref{(2.4)}.
\endproof

Taking into consideration the decomposition of $TM^{\bot}$, it can be easily proved:

\begin{proposition}
\label{split of fN}
Let $M$ be a semi-invariant
$\xi^{\bot}$-submanifold of a G.Q.S manifold $\widetilde{M}$. Then for any $N\in \Gamma(TM^{\bot})$ one has:
\begin{itemize}
\item[\emph{a)}] $BN\in\mathcal{D}^{\bot}$,
\item[\emph{b)}] $CN\in\mu$.
\end{itemize}
\end{proposition}

\begin{proposition}
\label{p:2.2} If $M$ is a semi-invariant
$\xi^{\bot}$-submanifold of a G.Q.S manifold $\widetilde{M}$ then
\begin{equation}
\label{(2.6)}
A_{\omega Z}W=A_{\omega W}Z
\end{equation}
for any $Z,W\in \Gamma ({\mathcal{D}}^{\bot})$.
\end{proposition}

The following two results give necessary and sufficient conditions for the integrability of the two distributions.

\begin{theorem}
\label{t:2.1}
Let $M$ be a semi-invariant $\xi^{\bot}$-submanifold of a G.Q.S manifold $\widetilde{M}$. Then the distribution
${\mathcal{D}}^{\bot} $ is integrable.
\end{theorem}
\proof
Let $Z,W\in \Gamma ({\mathcal{D}}^{\bot })$. Then from \eqref{(2.2)}, \eqref{(2.5)} and \eqref{(2.6)} we deduce that
\[
t[Z,W]=A_{\omega Z}W-A_{\omega W}Z=0.
\]
Hence the conclusion.
\endproof

\begin{theorem}
\label{t:2.2}
If $M$ is a semi-invariant $\xi^{\bot}$-submanifold
of a G.Q.S manifold $\widetilde{M}$ then the distribution ${\mathcal{D}}$ is integrable if
and only if
\begin{equation}
h(tX,Y)-h(X,tY)=({\mathcal  L}_{\xi }\tilde g)(X,\phi Y)~\xi,\quad X,Y\in \Gamma ({\mathcal{D}}).
\end{equation}
\end{theorem}
\proof
The statement yields directly from \eqref{(1.3)} and \eqref{(2.5)}
\[
\omega ([X,Y])=h(X,tY)-h(tX,Y)+({\mathcal L}_{\xi }\tilde g)(X,\phi Y)~\xi.
\]
\endproof

Notice that the two results above are analogue those obtained in the Kenmotsu case in \cite{9} and for the
cosymplectic case in \cite{s:hs}.
See also \cite{Mun05} when the submanifold is tangent to the structure vector field of the Sasakian manifold.

Moreover, from \eqref{(2.4)} we deduce

\begin{proposition}
\label{p:2.3}
Let $M$ be a $\xi ^{\bot }$-semi-invariant submanifold of a G.Q.S manifold $\widetilde{M}$. Then
\begin{equation}
\label{(2.8)}
A_{\xi }X=\alpha X,\quad \nabla _{X}^{\bot }\xi =-\beta X,\quad X\in \Gamma (TM).
\end{equation}
\end{proposition}

Let now $\{e_i, \phi e_i, e_{2p+j}\}$, $i\in\{1,...,p\}$, $j\in\{1,...,q\}$ be an adapted orthonormal local frame on $M$,
 where $q=\dim {\mathcal{D}}^{\bot}$ and $2p=\dim {\mathcal{D}}$. One can state the following
\begin{theorem}
\label{t:2.3}
If $M$ is a $\xi^{\bot}$-semi-invariant submanifold of a G.Q.S manifold $\widetilde{M}$ one has
\begin{equation*}
\label{(2.9)}
\eta (H)=\frac{1}{m}~{\rm trace}(A_{\xi }),\ \ m=2p+q.
\end{equation*}
\end{theorem}
\proof
Using a general formula for the mean curvature, e.g. $H=\frac{1}{m}\displaystyle\sum_{a=1}^{q}{\rm trace}\big(A_{\xi_a}\big)\xi_a$,
where $\{\xi_1,\ldots,\xi_q\}$ is an orthonormal basis in $TM^\bot$, the conclusion holds by straightforward computations.
\endproof

In the case when the ambient space is a Kenmotsu manifold we retrieve the known result from \cite[p. 614]{9}.
\begin{corollary}
\label{t:2.4}
There does not exist a minimal semi-invariant $\xi^{\bot}$-submanifold of a Kenmotsu manifold.
\end{corollary}

Also it is not difficult to prove:

\begin{theorem}
\label{t:2.5}
Let $M$ be a semi-invariant $\xi^{\bot}$-submanifold of a G.Q.S manifold $\widetilde{M}$. Then
\begin{enumerate}
\item the distribution ${\mathcal{D}}$ is integrable and its leaves are totally
geodesic in $M$ if and only if $h(X,Y)\in \Gamma (\mu )$, where $X,Y$ belong to ${\mathcal{D}}$;
\item any leaf of the integrable distribution ${\mathcal{D}}^{\bot }$ is totally geodesic in $M$
if and only if $h(X,Z)\in \Gamma (\mu )$ if $X\in
\Gamma ({\mathcal{D}})$ and $Z\in \Gamma ({\mathcal{D}}^{\bot })$.
\end{enumerate}
\end{theorem}
\proof
Let us prove only the first statement. For any $Z\in{\mathcal{D}}^\bot$ we have

\qquad
$\tilde g\big(h(X,Y),\phi Z\big)=\tilde g\big(\widetilde\nabla_XY,\phi Z\big)=-\tilde g\big(Y,\widetilde\nabla_X(\phi Z)\big)=$

\qquad\qquad\qquad
$=-\tilde g\big(Y,(\widetilde\nabla_X\phi)Z\big)-\tilde g\big(\phi Y,\widetilde\nabla_XZ\big)
=g\big(\nabla_X(\phi Y),Z\big)$.

Let $M^*$ be a leaf of the integrable distribution $\mathcal{D}$ and $h^*$ the second fundamental form of $M^*$ in $M$.

For any $Z\in \Gamma (\mathcal{D}^{\bot })$ we have:
$$
g(h^*(X, Y), Z)=\tilde{g}(\tilde{\nabla }_XtY, Z)=\tilde{g}((\tilde{\nabla }_X\varphi )Y+\varphi (\tilde{\nabla }_XY), Z)=-\tilde{g}(h(X, Y), \varphi Z)
$$
which proves that the leaf $M^*$ of the integrable $\mathcal{D}$ is totally geodesic in $M$ if and only if $h(X, Y)\in \Gamma (\mu )$.

\medskip

Notice that the part $(2)$ of the previous Theorem was obtained in the Kenmotsu case by Papaghiuc in \cite[p. 115]{p:iasi}.

\endproof

We end this section with the following

\begin{corollary}
\label{c:2.1}
If the leaves of the integrable distribution ${\mathcal{D}}$
are totally geodesic in $M$ then the structure vector field $\xi $
is ${\mathcal{D}}$-Killing, that is $({\mathcal  L}_{\xi }g)(X,Y)=0$, $X,Y\in \Gamma ({\mathcal{D}})$.
\end{corollary}

\section{Totally umbilical semi-invariant $\xi^{\bot}$-submanifolds}

The main purpose of this section is to obtain a complete characterization
of a totally umbilical semi-invariant $\xi^{\bot}$-submanifold of a G.Q.S
manifold $\widetilde{M}$. Recall that for a totally umbilical submanifold we have
\[
h(X, Y) = g(X, Y)H, \ \ X, Y\in\Gamma(TM).
\]

First we state:

\begin{theorem}
\label{t:3.1}
An invariant $\xi^{\bot}$-submanifold $M$ of a G.Q.S manifold is totally umbilical if and only if
\begin{equation}
\label{(3.1)}
h(X,Y)=\frac{1}{m}g(X,Y){\rm trace}\big(A_{\xi }\big)\xi .
\end{equation}
\end{theorem}
\proof
If $M$ is an invariant $\xi^{\bot}$-submanifold then for any $X,Y\in \Gamma (TM)$ we have
$h(X,\phi Y)=\phi h(X,Y)-g(A_{\xi }\phi X,Y)\xi $. Let us consider an orthonormal frame
\linebreak
$\{e_{i},e_{p+i}\}$, $i=1,\ldots,p$ on $M$; from the above relation one obtains
that $\phi H=0$. Again, since $M$ is an invariant submanifold:
\begin{equation}
\label{(3.2)}
H=g(H,\xi )\xi =\frac{1}{m}\sum_{i=1}^{m}g(h(e_{i},e_{i}),\xi )\xi=
   \frac{1}{m}{\rm trace}\big(A_{\xi }\big)\xi
\end{equation}
and the proof is complete.
\endproof
\begin{corollary}
\label{c:3.1}
A semi-invariant $\xi^{\bot}$-submanifold of a quasi-Sasakian manifold is minimal.
\end{corollary}

The case of a semi-invariant $\xi^{\bot}$-submanifold in a G.Q.S manifold
$\widetilde{M}$ is solved in the next Theorem.

\begin{theorem}
\label{t:3.2}
Let $M$ be a semi-invariant $\xi^{\bot}$-submanifold of a G.Q.S manifold $\widetilde{M}$ with $\dim {\mathcal{D}}^{\bot }> 1$.
Then $M$ is totally umbilical if and only if \eqref{(3.1)} holds.
\end{theorem}
\proof
Let $X\in \Gamma ({\mathcal{D}})$ be a unit vector field and $N\in \Gamma (\mu )\setminus{\rm span}\{\xi\}$. By direct
calculation it results that:
\[
g(H,N)=g(h(X,X),N)=g(\widetilde{\nabla}_{X}\phi X-(\widetilde{\nabla}_{X}\phi )X,\phi N)=g(h(X,\phi X),\phi N)=0
\]
which proves that $H\in \phi {\mathcal{D}}^{\bot }\oplus{\rm span}\{\xi\}$.

For $Z,W\in \Gamma ({\mathcal{D}}^{\bot })$, from \eqref{(2.5)} we derive $QA_{\phi Z}W=-g(Z,W)\phi H$ i.e.
\begin{equation}
\label{(3.3)}
g(Z,\phi H)g(W,\phi H)=g(Z,W)g(\phi H,\phi H).
\end{equation}
If we take $Z=W$ orthogonal to $\phi H$, since $\dim {\mathcal{D}}^{\bot }>1$, from
the above relation we infer $\phi H=0\Rightarrow H\in {\rm span}\{\xi \}$.
At this point the conclusion is straightforward.

Conversely, if \eqref{(3.1)} is supposed to be true, then we get \eqref{(3.2)} which together with \eqref{(3.1)} we deduce that
$M$ is totally umbilical.
\endproof

Let us remark that when $\widetilde{M}$ is a Kenmotsu manifold the
result of the Theorem~\ref{t:3.2} was proved in \cite{9}.

\begin{corollary}
\label{c:3.2}
Every $\xi ^{\bot }$-hypersurface of a G.Q.S manifold $\widetilde{M}$ is totally umbilical.
\end{corollary}
\proof
If $M$ is a hypersurface then $TM^{\bot }={\rm span}\{\xi \}$ that is
$h(X,Y)\in {\rm span}\{\xi \}$. Next, from \eqref{(3.2)} it follows \eqref{(3.1)}.
\endproof

In the particular case of a Kenmotsu manifold this result was obtained
 by Papaghiuc in \cite[p. 617]{9}.

\medskip

As a consequence of Theorem~\ref{t:3.2}, we obtain

\begin{theorem}
If $M$ is a totally umbilical semi-invariant
$\xi^{\bot}$-submanifold of a G.Q.S manifold $\widetilde{M}$ with $\dim\mathcal{D}^{\bot }> 1$,
then $M$ is a semi-invariant product.
\end{theorem}

Here, by a semi-invariant product we mean a semi-invariant $\xi^{\bot}$-submanifold of
$\widetilde{M}$ which can be locally written as a Riemannian product of a $\phi$-invariant
submanifold and a $\phi$-anti-invariant submanifold of $\widetilde{M}$, both of them orthogonal to $\xi$.

\proof
From the definition of totally umbilical submanifold we have $h(X,Z)=0$
for any $X\in\Gamma(\mathcal{D})$ and $Z\in\Gamma(\mathcal{D}^{\bot})$, so that, by b) of
Theorem~\ref{t:2.5}, the leaves of $\mathcal{D}^{\bot}$ are totally geodesic submanifolds of $M$.
By Theorem~\ref{t:3.2}, we have $h(X,Y)\in \mathrm{span}\{\xi\}\subset\mu$ for any $X,Y\in\mathcal  D$.
By virtue of a) of Theorem~\ref{t:2.4}, this implies that the invariant distribution $\mathcal  D$ is
integrable and its integral manifolds are totally geodesic submanifolds of $M$. Therefore, we conclude that $M$ is a semi-invariant product.
\endproof

Without any restriction on the dimension of $\mathcal{D}^{\bot}$, we have the following

\begin{theorem}
Let $M$ be a totally umbilical semi-invariant
$\xi^{\bot}$-submanifold of a G.Q.S manifold $\widetilde{M}$.
If $\mathcal  D$ is integrable, then each leaf of $\mathcal  D$ is a totally geodesic submanifold of $M$.
\end{theorem}

\proof
By using b) of Proposition~\ref{p:2.1}, for any $X\in\Gamma(\mathcal  D)$, we have
$$
\omega(\nabla _{X}X) = - g(X,X)CH - g(FX,\phi Y)\xi.
$$
Since $CH\in\mu$ by b) of Lemma \ref{split of fN} and $\omega U\in\phi\mathcal{D}^{\bot}$ for any
$U\in\Gamma(TM)$, from the above equation we deduce that $\omega(\nabla _{X}X)=0$, or equivalently
$$
\nabla _{X}X\in\mathcal{D}, \quad \forall X\in\Gamma(\mathcal{D}).
$$
Replacing $X$ by $X+Y$, we get $\nabla _{X}Y+\nabla _{Y}X\in\Gamma(\mathcal  D)$ for all $X,Y\in\Gamma(\mathcal  D)$.
This condition, together with the integrability of $\mathcal  D$, implies
\begin{equation}
\label{TG}
\nabla _{X}Y\in \mathcal{D}, \quad \forall X,Y\in\Gamma(\mathcal{D}).
\end{equation}
As $\mathcal  D$ is integrable, Frobenius theorem ensures that $M$ is foliated by leaves of $\mathcal  D$.
Combining this fact with (\ref{TG}), we conclude that the leaves of $\mathcal  D$ are totally geodesic submanifolds of $M$.
\endproof

{\small

\end{document}